\documentclass[12pt,reqno]{amsart}

\usepackage{amssymb}
\usepackage{amsfonts}

\newtheorem{coro}{Corollary}[section]
\newtheorem{exam}{Example}[section]
\newtheorem{theo}{Theorem}[section]
\newtheorem{lemm}{Lemma}[section]

\numberwithin{equation}{section}

\begin{document}
\begin{center}
{\large\bf Annuities under random rates of interest -- revisited}
\end{center}
\begin{center}
\vspace*{0.8cm} {Krzysztof BURNECKI \footnote{Research supported by KBN Grant No. PBZ KBN 016/P03/1999.},
Agnieszka MARCINIUK \\and Aleksander WERON\footnote{Corresponding author. Tel.: +48 71 3203530; fax: +48 71 3202654.\\
\indent {\it E-mail} address: weron@im.pwr.wroc.pl (A. Weron).}}\vskip 0.5cm {\it Hugo Steinhaus Center for
Stochastic Methods,\\Institute of Mathematics, Wroc{\l}aw University of Technology, Poland}
\end{center}
\vskip 1cm \centerline{\bf Abstract} \vskip 0.3cm In the article we consider accumulated values of
annuities-certain with yearly payments with independent random interest rates. We focus on annuities with
payments varying in arithmetic and geometric progression which are important basic varying annuities (see
Kellison, 1991). They appear to be a generalization of the types studied recently by Zaks (2001). We derive, via
recursive relationships, mean and variance formulae of the final values of the annuities. As a consequence, we
obtain moments related to the already discussed cases, which leads to a correction of main results from Zaks
(2001). \vskip 0.5cm \noindent {\it Keywords:} Finance mathematics; Annuity; Accumulated value; Random interest
rate

\section{Introduction}
An annuity is defined as a sequence of payments of a limited duration which we denote by $n$ (see, e.g. Gerber,
1997). The accumulated or final values of annuities are of our interest. Typically, for simplicity, it is assumed
that underlying interest rate is fixed and the same for all years. However, the interest rate that will apply in
future years is of course neither known nor constant. Thus, it seems reasonable to let interest rates vary in a
random way over time, cf. e.g. Kellison (1991).

We assume that annual rates of interest are independent random variables with common mean and variance.  We apply
this assumption in order to compute, via recursive relationships, fundamental characteristics, namely mean and
variance, of the accumulated values of annuities with payments varying in arithmetic and geometric progression
(see, e.g. Kellison, 1991). Since these important basic varying annuities can be reduced to the cases considered
by Zaks (2001), as a by-product we discover several mistakes in main results of Zaks (2001). Of course, {\em
errare humanum est}, but for the benefit of the readers we correct all of them here.

In Section 2 we introduce basic principles of the theory of annuities. Under the fixed interest assumption we
consider accumulated values of standard and non-standard annuities. We analyze different payment scenarios which
leads to contradiction to Corollary 2.1 from Zaks (2001). Finally, we introduce the ones with payments varying
according to arithmetic and geometric progression. It appears that all important types of annuities (cf.
Kellison, 1991) can be obtained as examples of the introduced ones.

In Section 3 we drop the assumption of fixed interest rates and we study the final values of the varying payment
annuities under stochastic approach to interest. We consider annual rates of interest to be independent random
variables with common mean and variance. Using recursive relations we compute first and second moment as well as
variance of the accumulated values. Special cases of the derived results correct the corresponding ones from Zaks
(2001), namely Corollary \ref{coro1}  Theorem 4.2, Corollary \ref{coro2} Theorem 4.3, Corollary \ref{coro5}
Theorem 4.5  and Corollary \ref{coro4} Theorem 4.6, respectively.

\section{Annuities with a fixed interest rate}

 First let us recall basic notation used in the theory of annuities.
Suppose that $j$ is the annual interest rate and fixed through the period of $n$ years. The annual discount rate
$d$ is given by the formula
\begin{equation}
(1+j)d=j \label{equat1}
\end{equation}
and the annual discount factor $v$ is given by the equation
\begin{equation}
(1+j)v=1. \label{equat2}
\end{equation}
Hence we have that
\begin{equation}
v+d=1. \label{equat3}
\end{equation}

In the article we concentrate on final or accumulated values of annuities. We assume that $k\leq{n}$ throughout,
unless otherwise specified. The accumulated value of an annuity-due with $k$ annual payments of $1$ is denoted by
$\ddot{s}_{\overline{k}|j}$ and given by the formulae
\begin{equation}
\ddot{s}_{\overline{k}|j}=\frac{(1+j)^k-1}{d} \label{equat4}
\end{equation}
and
\begin{equation}
\ddot{s}_{\overline{k}|j} =(1+j)^k+(1+j)^{k-1}+\ldots+(1+j)=(1+j)(1+ \ddot{s}_{\overline{k-1}|j}), \label{equat5}
\end{equation}
where the latter defines the recursive equation for $\ddot{s}_{\overline{k}|j}$.

Let us now consider a standard increasing annuity-due. The accumulated value of such annuity with $k$ annual
payments of $1,2,\ldots,k$, respectively is:
\begin{equation}
 (I{\ddot{s}})_{\overline{k}|j}=\frac{\ddot{s}_{\overline{k}|j}-k}{d}.
\label{equat6}
\end{equation}
This can be expressed recursively as
\begin{eqnarray}
(I{\ddot{s}})_{\overline{k}|j}&=&(1+j)^k+2(1+j)^{k-1}+\ldots+k(1+j)\nonumber\\
&=&(1+j)(k+(I{\ddot{s}})_{\overline{k-1}|j}), \label{equat7}
\end{eqnarray}
which corrects formula (2.8) from Zaks (2001).

The accumulated value of an increasing annuity-due with $k$ annual payments of $1^2, 2^2, \ldots, k^2$,
respectively is denoted by $({I^2}{\ddot{s}})_{\overline{k}|j}$ and calculated:
\begin{equation}
({I^2}{\ddot{s}})_{\overline{k}|j}=\frac{2(I{\ddot{s}})_{\overline{k}|j}- \ddot{s}_{\overline{k}|j}-k^2}{d}.
\label{equat8}
\end{equation}

The following two equations give the recursive formula for $({I^2}{\ddot{s}})_{\overline{k}|j}$ and set the
relationship between $({I^2}{\ddot{s}})_{\overline{k}|j}$ and $\ddot{s}_{\overline{k}|j}$.

\begin{eqnarray}
({I^2}{\ddot{s}})_{\overline{k}|j}&=&(1+j)^k+{2^2}(1+j)^{k-1} +\ldots+{k^2}(1+j)\nonumber\\
   &=&(1+j)(k^2+ ({I^2}{\ddot{s}})_{\overline{k-1}|j}),
\label{equat9}
\end{eqnarray}
\begin{equation}
({I^2}{\ddot{s}})_{\overline{k}|j}=\frac{(1+v)(\ddot{s}_{\overline{k}|j}+k^2)-2k- 2k^2}{d^2}. \label{equat10}
\end{equation}
The latter contradicts Corollary 2.1 from Zaks (2001).

In the sequel we will need the following two relations:
\begin{equation}
(I{\ddot{s}})_{\overline{k-1}|j}=(I{\ddot{s}})_{\overline{k}|j}-{\ddot{s}}_{\overline{k}|j}, \label{equat11}
\end{equation}

\begin{equation}
({I^2}{\ddot{s}})_{\overline{k-1}|j}=({I^2}{\ddot{s}})_{\overline{k}|j}- 2(I{\ddot{s}})_{\overline{k}|j}
+{\ddot{s}}_{\overline{k}|j}. \label{equat12}
\end{equation}

Standard decreasing annuities are similar to increasing ones, but the payments are made in the reverse order. The
accumulated value of such annuity-due with $k$ annual payments of $n, n-1, \ldots, n-k+1$, respectively is
denoted by $(D{\ddot{s}})_{\overline{n,k}|j}$ and given by the formulae:
\begin{eqnarray}
(D{\ddot{s}})_{\overline{n,k}|j}&=&n(1+j)^k+(n-1)(1+j)^{k-1}+\ldots+(n-k+1)(1+j)\nonumber\\
&=& (1+j)((D{\ddot{s}})_{\overline{n,k-1}|j}+(n-k+1)) \label{equat13}
\end{eqnarray}
and
\begin{equation}
(D{\ddot{s}})_{\overline{n,k}|j}=(n+1){\ddot{s}_{\overline{k}|j}}-(I{\ddot{s}})_{\overline{k}|j}. \label{equat14}
\end{equation}
The sum of a standard increasing annuity and its corresponding standard decreasing annuity is of course a
constant annuity.

Now let us consider the accumulated value of an annuity-due with payments varying in arithmetic progression (see,
e.g. Kellison, 1991). The first payment is $p$ and they subsequently increase by $q$ per period, i.e. they form a
sequence $(p, p+q, p+2q, \ldots, p+(k-1)q)$. We note that $p$ must be positive but $q$ can be either positive or
negative as long as $p+(k-1)q>0$ in order to avoid negative payments. The accumulated value of such annuity will
be denoted by $({\ddot{s}}_a)_{\overline{k}|j}^{(p,q)}$ and is defined by
\begin{equation}
({\ddot{s}}_a)_{\overline{k}|j}^{(p,q)}=p(1+j)^k+(p+q)(1+j)^{k-1}+\ldots+(p+(k- 1)q)(1+j). \label{equat15}
\end{equation}

Simple calculations lead to the following relationship.
\begin{equation}
({\ddot{s}}_a)_{\overline{k}|j}^{(p,q)}=(p- q){\ddot{s}}_{\overline{k}|j}+q(I{\ddot{s}})_{\overline{k}|j}.
\label{equat16}
\end{equation}

Important special cases are the combinations of $p=1$ and $q=0$, $p=1$ and $q=1$, and $p=n$ and $q=-1$.
\begin{exam}
\label{exam1} If $p=1$ and $q=0$, then $({\ddot{s}}_a)_{\overline{k}|j}^{(p,q)}$ becomes the accumulated value of
an annuity-due with $k$ annual payments of $1$, namely
\begin{equation}
({\ddot{s}}_a)_{\overline{k}|j}^{(1,0)}={\ddot{s}}_{\overline{k}|j}\label{equat17}.
\end{equation}
\end{exam}
\begin{exam}
\label{exam2} If $p=1$ and $q=1$, then $({\ddot{s}}_a)_{\overline{k}|j}^{(p,q)}$ becomes the accumulated value of
an increasing annuity-due with $k$ annual payments of $1, 2, \ldots, k$, respectively, namely
\begin{equation}
({\ddot{s}}_a)_{\overline{k}|j}^{(1,1)}=(I{\ddot{s}})_{\overline{k}|j}\label{equat18}.
\end{equation}
\end{exam}
\begin{exam}
\label{exam3} If $p=n$ and $q=-1$, then $({\ddot{s}}_a)_{\overline{k}|j}^{(p,q)}$ becomes the accumulated value
of a decreasing annuity-due with $k$ annual payments of $n, n-1, \ldots, n-k+1$, respectively, namely
\begin{equation}
({\ddot{s}}_a)_{\overline{k}|j}^{(n,-1)}=(D{\ddot{s}})_{\overline{n,k}|j}\label{equat19}.
\end{equation}
\end{exam}

Let us finally consider the accumulated value of an annuity-due with $k$ annual payments varying in geometric
progression. The first payment is $p$ and they subsequently increase in geometric progression with common ratio
$q$ per period, i.e. they form a sequence $(p, pq^2, pq^3, \ldots, pq^{k-1})$. We note that $p$ and $q$  must be
positive in order to avoid negative payments. The accumulated value of such annuity will be denoted by
$(\ddot{s}_g)_{\overline{k}|j}^{(p,q)}$ and is expressed as
\begin{eqnarray}
(\ddot{s}_g)_{\overline{k}|j}^{(p,q)}&=&p(1+j)^k+pq(1+j)^{k-1}+pq^2(1+j)^{k-
2}+\ldots+pq^{k-1}(1+j)\nonumber\\
&=&p(1+j){\frac{(1+j)^k-q^k}{1+j-q}}. \label{equat21}
\end{eqnarray}

Important special cases are the combinations of $p=1$ and $q=1$ (cf. Example \ref{exam1}), and $p=1$ and $q=1+u$.
\begin{exam}
\label{exam4} If $p=1$ and $q=1$, then $(\ddot{s}_g)_{\overline{k}|j}^{(p,q)}$ becomes the accumulated value of
an annuity-due with $k$ annual payments of 1, namely
\begin{equation}
({\ddot{s}}_g)_{\overline{k}|j}^{(1,1)}= {\ddot{s}}_{\overline{k}|j}\label{equat23}.
\end{equation}
\end{exam}

\begin{exam}
\label{exam5} If $p=1$ and $q=1+u$, then $(\ddot{s}_g)_{\overline{k}|j}^{(p,q)}$ becomes the accumulated value of
an annuity-due with $k$ annual payments of $1,1+u,(1+u)^2,\ldots,(1+u)^{k-1}$, respectively, namely
\begin{equation}
({\ddot{s}}_g)_{\overline{k}|j}^{(1,1+u)}= {\frac{(1+j)^k}{1+t}}{\ddot{s}}_{\overline{k}|t}\label{equat24},
\end{equation}
where $t$ is defined as the solution of
\begin{equation}
1+u=(1+j)(1+t)\label{equat25}.
\end{equation}
\end{exam}

\section{Annuities with random interest rates}

Let us suppose that the annual rate of interest in the $k$th year is a random variable $i_k$. We assume that, for
each $k$, we have $E(i_k)=j$ and $Var(i_k)=s^2$, and that $i_1, i_2, \ldots, i_n$ are independent random
variables. We write
\begin{equation}
E(1+i_k)=1+j=\mu \label{equat26}
\end{equation}
and
\begin{equation}
E[(1+i_k)^2]=(1+j)^2+s^2=1+f=m, \label{equat27}
\end{equation}
where
\begin{equation}
f=2j+j^2+s^2. \label{equat28}
\end{equation}
Obviously
\begin{equation}
Var(1+i_k)=m-{\mu}^2. \label{equat29}
\end{equation}

Next we define $r$ to be the solution of
\begin{equation}
1+r=\frac{1+f}{1+j} \label{equat30}
\end{equation}
and using (\ref{equat28}) we have
\begin{equation}
r=j+\frac{s^2}{1+j}. \label{equat31}
\end{equation}

For a $k$-year variable annuity-due with annual payments of $c_1, c_2, \ldots, c_k$, respectively, we denote
their final value by $C_k$.

\subsection{Payments varying in arithmetic progression}

In the case of payments varying in arithmetic progression we have that $c_k=p+\linebreak[3](k-1)q$, where $k=1,
2, \ldots, n$.

The final value of an annuity with such payments is given recursively:
\begin{equation}
C_k=(1+i_k)[C_{k-1}+(p+(k-1)q)]\;\;\;\;{\rm for}\;\;k=2,\ldots,n. \label{equat32}
\end{equation}

We can easily find $\mu_k=E(C_k)$ as
\begin{eqnarray}
E(C_k)&=&E\big((1+i_k)[C_{k-1}+(p+(k-1)q)]\big)\nonumber\\
&=&{E(1+i_k)}{E\big(C_{k-1}+(p+(k-1)q)\big)}
\end{eqnarray}
from independence of interest rates. Thus we have the recursive equation for $k=2, \ldots, n$:
\begin{equation}
\mu_k=\mu[{\mu}_{k-1}+(p+(k-1)q)]. \label{equat33}
\end{equation}

We note that $\mu_1=p(1+j)=p{\mu}$. The following lemma stems from (\ref{equat33}) and (\ref{equat15}).
\begin{lemm}
If $C_k$ denotes the final value of an annuity-due with annual payments varying in arithmetic progression: $p,
p+q, p+2q, \ldots, p+(k-1)q$, respectively and if the annual rate of interest during the $k$th year is a random
variable $i_k$ such that $E(1+i_k)=1+j$ and $Var(1+i_k)=s^2$, and $i_1, i_2, \ldots, i_n$ are independent random
variables, then
\begin{equation}
{\mu}_k=E(C_k)=({\ddot{s}}_a)_{\overline{k}|j}^{(p,q)}\label{equat34}.
\end{equation}
\label{tw1}
\end{lemm}

Similarly for the second moment $E(C_k^2)$ we have the recursive equation for $k=2, \ldots, n$:
\begin{equation}
m_k=E(C_k^2)=m[m_{k-1}+2(p+(k-1)q){\mu}_{k-1}+(p+(k-1)q)^2]. \label{equat333}
\end{equation}

We note that $m_1=p^2m$. In order to compute the second moment we need the following lemma.
\begin{lemm}
\label{lemm5} Under the assumptions of Lemma \ref{tw1} we have
\begin{equation}
m_k=M_{1k}+2M_{2k},
\end{equation}
where
\begin{equation}
M_{1k}= p^2m^k+(p+q)^2m^{k-1}+\ldots+(p+(k-1)q)^2m \label{equat36}
\end{equation}
and
\begin{equation}
M_{2k}=(p+q)m^{k-1}{\mu}_1+(p+2q)m^{k-2}{\mu}_2+\ldots+(p+(k-1)q)m{\mu}_{k-1}. \label{equat37}
\end{equation}
\end{lemm}
\begin{proof}
We proceed by induction.
When $k=2$, this follows on the equation (\ref{equat333}), since $\mu_1=p(1+j)=p{\mu}$ and $m_1=p^2m$. Assuming
our result is true for a given $k$ $(2\leq k\leq n-1)$, it stems from formula (\ref{equat333}) that it is also
true for $k+1$. This concludes the proof by induction.
\end{proof}

Since, by (\ref{equat26}), $1+f=m$ we can easily find that
\begin{equation}
M_{1k}=p^2{\ddot{s}}_{\overline{k}|f}+2pq(I{\ddot{s}})_{\overline{k- 1}|f}+q^2(I^2{\ddot{s}})_{\overline{k-1}|f}.
\label{equat38}
\end{equation}

Now we can apply (\ref{equat11}) and (\ref{equat12}) in order to derive equivalent expression for $M_{1k}$.
\begin{lemm}
\label{lemm3}
\begin{equation}
M_{1k}=(p-q)^2{\ddot{s}}_{\overline{k}|f}+2q(p-
q)(I{\ddot{s}})_{\overline{k}|f}+q^2(I^2{\ddot{s}})_{\overline{k}|f}\label{equat39}.
\end{equation}
\end{lemm}

Now we shall determine $M_{2k}$ using  (\ref{equat16}), (\ref{equat37}) and the fact that $1+f=m$. Writing
\begin{eqnarray}
M_{2k}&=&(p+q)(1+f)^{k-1}[(p-
q){\ddot{s}}_{\overline{1}|j}+q(I{\ddot{s}})_{\overline{1}|j}]\nonumber\\
&+&(p+2q)(1+f)^{k-2}[(p-q){\ddot{s}}_{\overline{2}|j}+q(I{\ddot{s}})_{\overline{2}|j}]+\ldots\nonumber\\
&+&(p+(k-1)q)(1+f) [(p-q){\ddot{s}}_{\overline{k-1}|j}+q(I{\ddot{s}})_{\overline{k-1}|j}]\nonumber\\
&=&\frac{d(p-q)+q}{d^2}\Big[((p+q)(1+f)^{k-1}(1+j)\nonumber\\
&+&(p+2q)(1+f)^{k-
2}(1+j)^2+\ldots\nonumber+(p+(k-1)q)(1+f)(1+j)^{k-1})\nonumber\\
&-&((p+q)(1+f)^{k-1}+(p+2q)(1+f)^{k-2}+\ldots\nonumber\\
&+&(p+(k-1)q)(1+f))\Big]-
\frac{q}{d}\Big[(p+q)(1+f)^{k-1}\nonumber\\
&+&2(p+2q)(1+f)^{k-2}+\ldots+ (k-1)(p+(k-1)q)(1+f)\Big]
\end{eqnarray}
and applying (\ref{equat30}) we obtain the following results.
\begin{lemm}
\label{lemm4} Under the assumptions of Lemma \ref{tw1} we have
\begin{eqnarray}
M_{2k}&=&\frac{1}{d^2}\Big[(p-q)(d(p-q)+q)(1+j)^k{\ddot{s}}_{\overline{k}|r}\nonumber\\
&+&q(d(p-q)+q)(1+j)^k(I{\ddot{s}})_{\overline{k}|r}\nonumber\\
&-&(p-q)(d(p-q)+qv)
{\ddot{s}}_{\overline{k}|f}\nonumber\\
&-&q(2d(p-q)+qv)(I{\ddot{s}})_{\overline{k}|f}- q^2d(I^2{\ddot{s}})_{\overline{k}|f}\Big]\label{equat40}.
\end{eqnarray}
\end{lemm}

\begin{lemm}
\label{lemm1} Under the assumptions of Lemma \ref{tw1} we have
\begin{eqnarray}
m_k&=&\frac{1}{d^2}\Big[(q-p)(d(p-q)(1+v)+2qv){\ddot{s}}_{\overline{k}|f}\nonumber\\
&-&2q(d(p-q)(1+v)+qv)(I{\ddot{s}})_{\overline{k}|f}\nonumber\\
&-&dq^2(1+v)(I^2
{\ddot{s}})_{\overline{k}|f}+2(p-q)(d(p-q)+q)(1+j)^k{\ddot{s}}_{\overline{k}|r}\nonumber\\
&+&2q(d(p-q)+q)(1+j)^k(I{\ddot{s}})_{\overline{k}|r}\Big]\label{equat41}.
\end{eqnarray}
\end{lemm}

We have thus reached a formula for $E(C_k^2)$. In order to compute $Var(C_k)$ we need yet an expression for
$E(C_k)^2$.

\begin{lemm} Under the assumptions of Lemma \ref{tw1} we have
\label{lemm2}
\begin{eqnarray}
{\mu}_k^2&=&\frac{p-q}{d}\Big(p-q+\frac{2q}{d}\Big)\Big({\ddot{s}}_{\overline{2k}|j}-
2{\ddot{s}}_{\overline{k}|j}\Big)-\frac{2q(p-q)k}{d}{\ddot{s}}_{\overline{k}|j}\nonumber\\
&+&\left(\frac{q}{d}\right)^2\Big((I{\ddot{s}})_{\overline{2k}|j}-2(1+kd)(I{\ddot{s}})_{\overline{k}|j}-
k^2\Big)\label{equat42}.
\end{eqnarray}
\end{lemm}
\begin{proof}
It is easy to show that
\begin{equation}
(\ddot{s}_{\overline{k}|j})^2=\frac{\ddot{s}_{\overline{2k}|j}-2{\ddot{s}}_{\overline{k}|j}}{d} \label{equat43}
\end{equation}
and
\begin{equation}
(I{\ddot{s}})_{\overline{k}|j}^2=\frac{(I{\ddot{s}})_{\overline{2k}|j}-
2(1+kd)(I{\ddot{s}})_{\overline{k}|j}-k^2}{d^2}, \label{equat44}
\end{equation}
cf. Lemma 3.3 and 4.3 from Zaks (2001). From (\ref{equat16}) we may write
\begin{eqnarray}
{\mu}_k^2&=&((p-q){\ddot{s}}_{\overline{k}|j}+q(I{\ddot{s}})_{\overline{k}|j})^2\nonumber\\
&=&(p-q)^2({\ddot{s}}_{\overline{k}|j})^2+2q(q-p){\ddot{s}}_{\overline{k}|j}
(I{\ddot{s}})_{\overline{k}|j}+q^2(I{\ddot{s}})_{\overline{k}|j}^2.
\end{eqnarray}
Substituting from (\ref{equat43}), (\ref{equat44}) and (\ref{equat6}) completes the proof.
\end{proof}

Now, we are allowed to state the following theorem.
\begin{theo}
\label{theo1} Under the assumptions of Lemma \ref{tw1} we have
    \begin{eqnarray}
    E(C_k)&=&({\ddot{s}}_a)_{\overline{k}|j}^{(p,q)},\\
    Var(C_k)&=&m_k-{\mu}_k^2,
    \end{eqnarray}
    where $m_k$ is given by Lemma \ref{lemm1} and ${\mu}_k^2$ by Lemma \ref{lemm2}.
\end{theo}

Let us now consider the situation when $p=1$ and $q=0$. We know, from Example \ref{exam1}, that it is the case of
an annuity-due with $k$ annual payments of $1$. Then we obtain the following corollary (cf. Theorem 3.2 from
Zaks, 2001).
\begin{coro}
\label{coro3} If $C_k$ denotes the  final value of an annuity-due with $k$ annual payments of $1$ and if the
annual rate of interest during the $k$th year is a random variable $i_k$ such that $E(1+i_k)=1+j$ and
$Var(1+i_k)=s^2$, and $i_1, i_2, \ldots, i_n$ are independent random variables, then
\begin{equation}E(C_k)={\ddot{s}}_{\overline{k}|j},\end{equation}
\begin{equation}Var(C_k)=\frac{2(1+j)^{k+1}{\ddot{s}}_{\overline{k}|r}-(2+j){\ddot{s}}_{\overline{k}|f}-
(1+j){\ddot{s}}_{\overline{2k}|j}+2(1+j){\ddot{s}}_{\overline{k}|j}}{j}.\end{equation}
\end{coro}

Another important case is the combination of $p=1$ and $q=1$, see Example \ref{exam2}. It is an annuity-due with
$k$ annual payments of $1, 2, \ldots, k$. The following corollary is a direct consequence of Lemma \ref{tw1},
\ref{lemm3}, \ref{lemm4} and \ref{lemm1}.

\begin{coro}
\label{coro1} If $C_k$ denotes the  final value of an increasing annuity-due with $k$ annual payments of of $1,
2, \ldots, k$, respectively and if the annual rate of interest during the $k$th  year is a random variable $i_k$
such that $E(1+i_k)=1+j$ and $Var(1+i_k)=s^2$, and $i_1, i_2, \ldots, i_n$ are independent random variables, then
\begin{enumerate}
\item[(a)] $E(C_k)=(I\ddot{s})_{\overline{k}|j}$,
\item[(b)] $M_{1k}=(I^2\ddot{s})_{\overline{k}|f}$,
\item[(c)] $M_{2k}=\frac{(1+j)^{k+2}(I\ddot{s})_{\overline{k}|r}-
(1+j)(I\ddot{s})_{\overline{k}|f}-j(1+j)(I^2\ddot{s})_{\overline{k}|f}}{j^2}$,
\item[(d)] $m_k=\frac{2(1+j)^{k+2}(I\ddot{s})_{\overline{k}|r}-
2(1+j)(I\ddot{s})_{\overline{k}|f}-j(2+j)(I^2\ddot{s})_{\overline{k}|f}}{j^2}$.
\end{enumerate}
\end{coro}
Part (b) of Corollary \ref{coro1} corrects (4.6) from Zaks (2001) and (d) is in contradiction to Theorem 4.2 from
Zaks (2001). These results can be summarized in the following corollary.
\begin{coro}
\label{coro2} Under the assumptions of Corollary \ref{coro1} we have
\begin{equation}E(C_k)=(I\ddot{s})_{\overline{k}|j},\end{equation}
\begin{eqnarray}
Var(C_k)&=&\frac{2(1+j)^{k+2}(I\ddot{s})_{\overline{k}|r}-
2(1+j)(I\ddot{s})_{\overline{k}|f}-j(2+j)(I^2\ddot{s})_{\overline{k}|f}}{j^2}\nonumber\\
&-&\frac{(I\ddot{s})_{\overline{2k}|j}-2(1+kd)(I\ddot{s})_{\overline{k}|j}-k^2}{d^2}.
\end{eqnarray}
\end{coro}
Corollary \ref{coro2} (variance part) contradicts Theorem 4.3 from Zaks (2001). Let us finally consider the
situation when $p=n$ and $q=-1$, see Example \ref{exam3}. Then we obtain the following corollary.
\begin{coro}
\label{coro5} If $C_k$ denotes the  final value of a decreasing annuity-due with $k$ annual payments of of $n,
n-1, \ldots, n-k+1$, respectively and if the annual rate of interest during the $k$th  year is a random variable
$i_k$ such that $E(1+i_k)=1+j$ and $Var(1+i_k)=s^2$, and $i_1, i_2, \ldots, i_n$ are independent random
variables, then
\begin{equation}E(C_k)=(D{\ddot{s}})_{\overline{n,k}|j},
\end{equation}
\begin{eqnarray}
Var(C_k)&=&\frac{\ell}{d^2} \Big[\frac{(n-1/j)^2(1+j)^{2k}\ddot{s}_{\overline{k}|\ell}}{1+\ell}
-\frac{2(n-1/j)^2 (1+j)^k\ddot{s}_{\overline{k}|r}}{1+r}\nonumber\\
&+&\frac{(n-1/j)^2\ddot{s}_{\overline{k}|f}}{1+f}+\frac{2(n-1/j)(1+j)^k
(I\ddot{s})_{\overline{k}|r}}{1+r}\nonumber\\
&-&\frac{2(n-1/j)(I\ddot{s})_{\overline{k}|f}}{1+f} +\frac{(I^2\ddot{s})_{\overline{k}|f}}{1+f}\Big],
\end{eqnarray}
where $\ell=(s/(1+j))^2$.
\end{coro}
Corollary \ref{coro5} (variance part) corrects Theorem 4.5 from Zaks (2001).

\subsection{Payments varying in geometric progression}

In the case of annuities-due with payments varying in geometric progression we have that $c_k=pq^{k-1}$, where
$k=1, 2, \ldots, n$.

The final value of that annuity is given recursively:
\begin{equation}
C_k=(1+i_k)[C_{k-1}+pq^{k-1}]\label{equat45}.
\end{equation}

Similarly, as in the case of payments varying in arithmetic progression, we easily find that for $k=2, \ldots, n$
\begin{equation}
\mu_k=E(C_k)={\mu}[{\mu}_{k-1}+pq^{k-1}]. \label{equat46}
\end{equation}

The second moment $E(C_k^2)$ is given by
\begin{equation}
m_k=E(C_k^2)=m[m_{k-1}+2pq^{k-1}{\mu}_{k-1}+p^2q^{2(k-1)}]. \label{equat47}
\end{equation}

We note that $\mu_1=p(1+j)=p{\mu}$ and $m_1=p^2m$. In analogy with Lemma \ref{tw1} we obtain a pleasing form of
$E(C_k)$.
\begin{lemm}
If $C_k$ denotes the  final value of an annuity-due with $k$ annual payments varying in geometric progression:
$p, pq, pq^2, \ldots, pq^{k-1}$, respectively and if the annual rate of interest during the $k$th year is a
random variable $i_k$ such that $E(1+i_k)=1+j$ and $Var(1+i_k)=s^2$, and $i_1, i_2, \ldots, i_n$ are independent
random variables, then
\begin{equation}
{\mu}_k=E(C_k)=({\ddot{s}}_g)_{\overline{k}|j}^{(p,q)}\label{equat48}.
\end{equation}
\label{tw2}
\end{lemm}

Similarly, as in the previous point, in order to find a formula for $Var(C_k)$ we are about to compute $m_k$ and
$\mu_k^2$. We commence by calculating $m_k$.
\begin{lemm}
Under the assumptions of Lemma \ref{tw2} we have
\begin{eqnarray}
m_k&=&p^2m^k+p^2q^2m^{k-1}+\ldots+p^2q^{2(k-1)}m\nonumber\\
&+&2[pqm^{k-1}{\mu}_1+pq^2m^{k- 2}{\mu}_2+\ldots+pq^{k-1}m{\mu}_{k-1}]\label{equat49}.
\end{eqnarray}
\end{lemm}
\begin{proof}
The thesis follows by induction, applying (\ref{equat47}) and the fact $1+f=m$.
\end{proof}

Let
\begin{equation}
M_{1k}= p^2m^k+p^2q^2m^{k-1}+\ldots+p^2q^{2(k-1)}m \label{equat50}
\end{equation}
and
\begin{equation}
M_{2k}=pqm^{k-1}{\mu}_1+pq^2m^{k-2}{\mu}_2+\ldots+pq^{k-1}m{\mu}_{k-1}. \label{equat51}
\end{equation}
Hence
\begin{equation}
m_k=M_{1k}+2M_{2k}
\end{equation}
(cf. Lemma \ref{lemm5}). Since $1+f=m$, we can easily obtain an elegant expression for $M_{1k}$.
\begin{lemm}
\label{lemm6}
\begin{equation}
M_{1k}=p^2(1+f)\frac{(1+f)^k-q^{2k}}{1+f- q^2}=({\ddot{s}_g})_{\overline{k}|f}^{(p^2,q^2)}\label{equat52}.
\end{equation}
\end{lemm}

Now we rewrite (\ref{equat51}) applying  $1+f=m$  and $1+f=(1+j)(1+r)$ giving
\begin{eqnarray}
M_{2k}&=&pq(1+f)^{k-1}p(1+j)\frac{(1+j)-q}{1+j-q}\nonumber\\
&+&pq^2(1+f)^{k- 2}p(1+j)\frac{(1+j)^2-q^2}{1+j-q}
+\ldots\nonumber\\
&+&pq^{(k-1)}(1+f) p(1+j)\frac{(1+j)^{k-
1}-q^{k-1}}{1+j-q}\nonumber\\
&=&\frac{p^2(1+j)}{1+j-q}\Big[(q(1+f)^{k-1}(1+j)+q^2(1+f)^{k-2}(1+j)^2+\dots\nonumber\\
&+&q^{k-
1}(1+f)(1+j)^{k-1}+(1+f)^k-(1+f)^k)-(q^2(1+f)^{k-1}\nonumber\\
&+&q^4(1+f)^{k-
2}+\ldots+q^{2(k-1)}(1+f)+(1+f)^k-(1+f)^k)\Big]\nonumber\\
&=&\frac{p^2(1+j)}{1+j-q}\left[(1+j)^k(1+r){\frac{(1+r)^k-q^k}{1+r-q}}-
\frac{({\ddot{s}}_q)_{\overline{k}|f}^{(p^2,q^2)}}{p^2}\right].
\end{eqnarray}

Therefore we may write the following lemma.
\begin{lemm}
\label{lemm7}
\begin{equation}
M_{2k}=\frac{p(1+j)^{k+1}({\ddot{s}}_g)_{\overline{k}|r}^{(p,q)}-
(1+j)({\ddot{s}}_g)_{\overline{k}|f}^{(p^2,q^2)}}{1+j-q}\label{equat53}
\end{equation}
\end{lemm}

By virtue of Lemma \ref{lemm6} and \ref{lemm7}, and the fact that $m_k=M_{1k}+2M_{2k}$ we have the following
lemma.
\begin{lemm}
Under the assumptions of Lemma \ref{tw2} we have
\begin{equation}
m_k=\frac{2p(1+j)^{k+1}({\ddot{s}}_g)_{\overline{k}|r}^{(p,q)}-
(q+1+j)({\ddot{s}}_g)_{\overline{k}|f}^{(p^2,q^2)}}{1+j-q}\label{equat54}.
\end{equation}
\end{lemm}

Thus, we have reached a formula for $E(C_k^2)$. Now we need to derive an expression for $E(C_k)^2$.
\begin{lemm}Under the assumptions of Lemma \ref{tw2} we have
\begin{equation}
{\mu}_k^2=\frac{p(1+j)}{1+j-q}\left(({\ddot{s}}_g)_{\overline{2k}|j}^{(p,q)}-
2q^k({\ddot{s}}_g)_{\overline{k}|j}^{(p,q)}\right)\label{equat55}.
\end{equation}
\end{lemm}
\begin{proof}
From Lemma \ref{tw2} and (\ref{equat21}), we have that
\begin{eqnarray}
{\mu}_k^2&=&\left(p(1+j)\frac{(1+j)^k-q^k}{1+j-q}\right)^2\nonumber\\
&=&p^2(1+j)^2\frac{(1+j)^{2k}-
2(1+j)^kq^k+q^{2k}}{(1+j-q)^2}\nonumber\\
&=& \frac{p^2(1+j)^2}{1+j-q}\left(\frac{(1+j)^{2k}-
q^{2k}}{1+j-q}-\frac{2q^k((1+j)^k-q^k)}{1+j-q}\right)\hspace{-0.2cm},
\end{eqnarray}
which using (\ref{equat21}) completes the proof.
\end{proof}

Since $Var(C_k)=m_k-{\mu}_k^2$, we have following theorem.
\begin{theo}
\label{theo2} Under the assumptions of Lemma \ref{tw2} we have
\begin{eqnarray}
Var(C_k)&=&\frac{2p(1+j)^{k+1}({\ddot{s}}_g)_{\overline{k}|r}^{(p,q)}-
(1+j+q)({\ddot{s}}_g)_{\overline{k}|f}^{(p^2,q^2)}}
{1+j-q}\nonumber\\
&-&\frac{p(1+j)\left(({\ddot{s}}_g)_{\overline{2k}|j}^{(p,q)}-
2q^k({\ddot{s}}_g)_{\overline{k}|j}^{(p,q)}\right)}{1+j-q}.
\end{eqnarray}
\end{theo}

An important case, see Example \ref{exam4}, is the combination of $p=1$ and $q=1$. Then we obtain an annuity-due
with $k$ annual payments of $1$ and Theorem \ref{theo2} yields Corollary \ref{coro3}.

Another important case is the combination of $p=1$ and $q=1+u$, where $1+u=(1+j)(1+t)$. This defines an
annuity-due with $k$ annual payments of $1, 1+u, (1+u)^2, \ldots, (1+u)^{k-1}$, respectively. We assume also that
$1+f=(1+u)^2(1+h)$ and $1+f=(1+j)^2(1+t)(1+w)$. This leads to the following corollary.
\begin{coro}
\label{coro4} If $C_k$ denotes the  final value of an annuity-due with $k$ annual payments of $1, 1+u, (1+u)^2,
\ldots, (1+u)^{k-1}$, respectively and if the annual rate of interest during the $k$th year is a random variable
$i_k$ such that $E(1+i_k)=1+j$ and $Var(1+i_k)=s^2$, and $i_1, i_2, \ldots, i_n$ are independent random
variables, then
$$E(C_k)=\frac{(1+j)^k{\ddot{s}}_{\overline{k}|t}}{1+t},$$
\begin{eqnarray*}
Var(C_k)&=&\frac{(1+u)^{2k}(2+t){\ddot{s}}_{\overline{k}|h}-
2(1+j)^{2k}(1+t)^k{\ddot{s}}_{\overline{k}|w}}{t}\\
&-&\frac{(1+j)^{2k}({\ddot{s}}_{\overline{2k}|t}-2{\ddot{s}}_{\overline{k}|t})}{t(1+t)}.
\end{eqnarray*}
\end{coro}

Corollary \ref{coro4} (variance part) contradicts Theorem 4.6 from Zaks (2001).

\end{document}